\documentclass[12pt]{article}

\usepackage{comment}
% ...

\newtheorem{thm}{Theorem}
\newtheorem{clm}[thm]{Claim}
\newtheorem{cor}[thm]{Corollary}

\newtheorem{lem}[thm]{Lemma}

\usepackage{amsmath,amsfonts}
\usepackage{graphicx}
\usepackage{tikz}
\usepackage{verbatim}

\def\f{{\phi}}

\def\s{{\sigma}}

\def\th{{\theta}}
\def\thp{{\th^\pr}}

\def\zN{{\mathbb N}}

\def\cP{{\cal P}}

\def\sNP{{\sf NP}}
\def\sP{{\sf P}}

\def\bp{{b^\pr}}

\def\Ch{{\hat{C}}}
\def\Cp{{C^\pr}}
\def\Cpp{{C^{\pr\pr}}}

\def\dH{{d_H}}
\def\dpr{{d^\pr}}
\def\dppr{{d^{\pr\pr}}}

\def\Gh{{\hat{G}}}
\def\Gp{{G^\pr}}
\def\Gpp{{G^{\pr\pr}}}
\def\Gpi{{G^\pr_i}}
\def\Gpj{{G^\pr_j}}
\def\Gpo{{G^\pr_1}}
\def\Gpt{{G^\pr_2}}

\def\nH{{n_H}}
\def\np{{n^\pr}}
\def\npp{{n^{\pr\pr}}}

\def\Nx!{{N_x^!}}

\def\rp{{r^\pr}}

\def\Tp{{T^\pr}}

\def\cost{{\sf cost}}

\def\diam{{\sf diam}}
\def\dist{{\sf dist}}

\def\ecc{{\sf ecc}}

\def\len{{\sf length}}

\def\mt{{\emptyset}}

\def\ppr{{\pr\pr}}
\def\pr{{\prime}}
\def\Pf{{\noindent\bf Proof.\ \ }}

\def\pf{{\hfill $\Box$\bigskip}}

\def\rar{{\rightarrow}}

\def\sse{{\subseteq}}

\begin{document}

\title{Pebbling in Semi-2-Trees}

\author{
Liliana Alc\' on\thanks{
Departamento de Matem\' atica,
FCE, UNLP, La Plata, Argentina
}\\
Marisa Gutierrez\footnotemark[1]
\ \thanks{
CONICET, Argentina
}\\
Glenn Hurlbert\thanks{
Department of Mathematics and Applied Mathematics,
Virginia Commonwealth University, USA
}
\thanks{
Research partially supported by Simons Foundation Grant \#246436.
}
}
\maketitle

\newpage

%
%
\begin{abstract}

Graph pebbling is a network model for transporting discrete resources
that are consumed in transit.  Deciding whether a given configuration
on a particular graph can reach a specified target is {\sf NP}-complete,
even for diameter two graphs, and deciding whether the pebbling
number has a prescribed upper bound is $\Pi_2^{\sf P}$-complete.
Recently we proved that the pebbling number of a split graph can
be computed in polynomial time.  This paper advances the program
of finding other polynomial classes, moving away from the large
tree width, small diameter case (such as split graphs) to small
tree width, large diameter, continuing an investigation on the
important subfamily of chordal graphs called $k$-trees.  In particular,
we provide a formula, that can be calculated in polynomial time,
for the pebbling number of any semi-2-tree, falling shy of the result
for the full class of 2-trees.

\smallskip

2xik
\end{abstract}

\newpage

%
%
\section{Introduction}\label{Intro}

The fundamental question in graph pebbling is whether a given supply
({\it configuration}) of discrete pebbles on the vertices of a
connected graph can satisfy a particular set of demands on the
vertices.  The operation of pebble movement across an edge $\{u,v\}$
is called a {\it pebbling step}: while two pebbles cross the edge,
only one arrives at the opposite end, as the other is consumed.  We
write $(u,v)$ to denote a pebbling step from $u$ to $v$.  The most
studied scenario involves the demand of one pebble on a single {\it
root} vertex $r$.  Satisfying this demand is often referred to as
{\it reaching} or {\it solving} $r$, and configurations are
consequently called either $r$-{\it solvable} or $r$-{\it unsolvable}.

The {\it size} $|C|$ of a configuration $C:V\rar\zN=\{0,1,\ldots\}$
is its total number of pebbles $\sum_{v\in V}C(v)$.  The {\it
pebbling number} $\pi(G)= \max_{r\in V}\pi(G,r)$, where $\pi(G,r)$
is defined to be the minimum number $s$ so that every configuration
of size at least $s$ is $r$-solvable.  Simple sharp lower bounds
like $\pi(G)\ge n$ and $\pi(G)\ge 2^{\diam(G)}$ are easily derived.
Graphs satisfying $\pi(G)=n$ are called {\it Class $0$} and are a
topic of much interest.  Recent chapters in \cite{HurlHGT} and
\cite{HurlMMC} include variations on the theme such as $k$-pebbling,
fractional pebbling, optimal pebbling, cover pebbling, and pebbling
thresholds, as well as applications to combinatorial number theory,
combinatorial group theory, and $p$-adic diophantine equations, and
also contain important open problems in the field.

Computing the pebbling number is difficult in general.  The
problem of deciding if a given configuration on a graph can reach
a particular vertex was shown in \cite{HurKie} and \cite{MilCla}
to be \sNP-complete, even for diameter two graphs (\cite{CuLeSiTa})
or planar graphs (\cite{LeCuDi}).  Interestingly, the problem was
shown in \cite{LeCuDi} to be in \sP\ for graphs that are both planar
and diameter two, as well as for outerplanar graphs (which include
2-trees).  The problem of deciding whether a graph $G$ has pebbling
number at most $k$ was shown in \cite{MilCla} to be $\Pi_2^{\sf
P}$-complete.

In contrast, the pebbling number is known for many graphs.  For
example, in \cite{PaSnVo} the pebbling number of a diameter $2$
graph $G$ was determined to be $n$ or $n+1$.  Moreover, \cite{ClHoHu}
and \cite{BlaSch} characterized those graphs having $\pi(G)=n+1$,
and it was shown in \cite{HeHeHu} that one can recognize such graphs
in quartic time, improving on the order $n^3m$ algorithm of \cite{BekCus}.
Beginning a program to study for which graphs
their pebbling number can be computed in polynomial time, the
authors of \cite{AlGuHu} produced a formula for the family of split
graphs that involves several cases.  For a given graph, finding to
which case it belongs takes $O(n^{1.41})$ time.  The authors also
conjectured that the pebbling number of a chordal graph of bounded
diameter can be computed in polynomial time.

In opposition to the small diameter, large tree width case of split
graphs, we turn here to chordal graphs with large diameter and small
tree width\footnote{One can find the definition of tree-width in 
\cite{B}, but it is not necessary for this paper.}.
In this paper we study 2-paths, the sub-class of 2-trees
whose graphs have exactly two simplicial vertices, as well as what
we call semi-2-trees, the sub-class of 2-trees, each of whose
blocks are 2-paths, and prove an exact formula that can be computed 
in linear time.

%
%

\section{Preliminary Definitions and Results}\label{Prelim}

In order to simplify notation, for a subgraph $H\subset G$ or
subset $H\subset V(G)$ we write $C(H)$ to denote 
$\sum_{v\in V(H)}C(v)$.
We use $C_H$ for the restriction of $C$ to $H$.

A {\it simplicial} vertex in a graph is a vertex whose neighbors
form a complete graph.  It is $k$-{\it simplicial} if it also has
degree $k$.  A $k$-{\it tree} is a graph $G$ that is either a
complete graph of size $k$ or has a $k$-simplicial vertex $v$ for
which $G-v$ is a $k$-tree.  A $k$-{\it path} is a $k$-tree with
exactly two simplicial vertices.  A {\it semi-$2$-tree} is a graph
in which each of its blocks is a 2-path, with each of its cut-vertices
being simplicial in all of its blocks.  For the purpose of our
work we derive a new characterization of 2-paths that facilitates
the analysis of its pebbling number.

Let $P=x_0,x_1,\ldots,x_{d-1},x_d$ be a shortest $rs$-path between
two vertices $r=x_0$ and $s=x_d$ of $G$, where
$d=\dist(r,s)=\diam(G)$. For $1\leq i \leq d-1$, an
$x_{i-1}x_{i+1}$-{\it fan} ({\it centered on} $x_i$) is a subgraph $F$ of $G$
consisting of the subpath $x_{i-1},x_i,x_{i+1}$ of $P$ and a path
$Q=x_{i-1},v_{i,1},\ldots,v_{i,k_i},x_{i+1}$ with $k_i\geq 1$ such
that $x_i$ is adjacent to every vertex of $Q$. We call $F^\pr$
the set $\{v_{i,1},\ldots,v_{i,k_i}\}$.

Let $F_i$  be an $x_{i-1}x_{i+1}$-fan and  $F_{i+1}$ be  an
$x_{i}x_{i+2}$-fan, centered on $x_i$ and on $x_{i+1}$,
respectively. We say that $F_i$ and $F_{i+1}$ are
\textit{opposite-sided} if $F_i^\pr \cap F_{i+1}^\pr=\emptyset$;
and that they are \textit{same-sided} when $F_i^\pr \cap
F_{i+1}^\pr=\{v_{i,k_i}\}$ and $v_{i,k_i}=v_{i+1,1}$.

\begin{figure}
\begin{center}

\begin{tikzpicture}[scale=.9]
\tikzstyle{every node}=[draw,circle,fill=black,minimum size=2pt,inner sep=3pt]
\draw (-5.0,0.0) node (r) [label=below: $r$] {};
\path (-4.0,0.0) node (x1) [label=below: $x_1$] {};
\path (-4.7,0.7) node (v11) [label=above: $v_{1,1}$] {};
\path (-4.0,1.0) node (v12) [label=above: $v_{1,2}$] {};
\path (-3.3,0.7) node (v13) [label={[xshift=0cm, yshift=-.2cm]$v_{1,3}$}] {};
\path (-3.0,0.0) node (x2) [label=below: $x_2$] {};
\path (-3.5,-.85) node (v21) [label=below: $v_{2,1}$] {};
\path (-2.5,-.85) node (v22) [label=below: $v_{2,2}$] {};
\path (-2.0,0.0) node (x3) [label=below: $x_3$] {};
\path (-2.7,0.7) node (v31) [label=above: $v_{3,1}$] {};
\path (-2.0,1.0) node (v32) [label=above: $v_{3,2}$] {};
\path (-1.3,0.7) node (v33) [label={[xshift=0cm, yshift=-.2cm]$v_{3,3}$}] {};
\path (-1.0,0.0) node (s) [label=below: $s$] {};
\draw (r)
  -- (x1)
  -- (x2)
  -- (x3)
  -- (s);
\draw (r)
  -- (v11)
  -- (v12)
  -- (v13)
  -- (x2);
\draw (x1)
  -- (v21)
  -- (v22)
  -- (x3);
\draw (x2)
  -- (v31)
  -- (v32)
  -- (v33)
  -- (s);
\draw (x1) -- (v11);
\draw (x1) -- (v12);
\draw (x1) -- (v13);
\draw (x2) -- (v21);
\draw (x2) -- (v22);
\draw (x3) -- (v31);
\draw (x3) -- (v32);
\draw (x3) -- (v33);

\draw (0.0,0.0) node (r) [label=below: $r$] {};
\path (1.0,0.0) node (x1) [label=below: $x_1$] {};
\path (0.3,0.7) node (v11) [label=above: $v_{1,1}$] {};
\path (1.0,1.0) node (v12) [label=above: $v_{1,2}$] {};
\path (2.0,0.7) node (v13) [label=above: $v_{1,3}$] {};
\path (2.0,0.0) node (x2) [label=below: $x_2$] {};
\path (2.0,0.7) node (v21) {};
\path (3.0,0.0) node (x3) [label=below: $x_3$] {};
\path (2.0,0.7) node (v31) {};
\path (3.5,0.85) node (v32) [label=above: $v_{3,2}$] {};
\path (4.0,0.0) node (s) [label=below: $s$] {};
\draw (r)
  -- (x1)
  -- (x2)
  -- (x3)
  -- (s);
\draw (r)
  -- (v11)
  -- (v12)
  -- (v13)
  -- (x2);
\draw (x1)
  -- (v21)
  -- (x3);
\draw (x2)
  -- (v31)
  -- (v32)
  -- (s);
\draw (x1) -- (v11);
\draw (x1) -- (v12);
\draw (x1) -- (v13);
\draw (x2) -- (v21);
\draw (x3) -- (v31);
\draw (x3) -- (v32);

\draw (5.0,0.0) node (r) [label=below: $r$] {};
\path (6.0,0.0) node (x1) [label=below: $x_1$] {};
\path (5.5,0.85) node (v11) [label=above: $v_{1,1}$] {};
\path (6.5,0.85) node (v12) [label=above: $v_{1,2}$] {};
\path (7.0,0.0) node (x2) [label=above: $x_2$] {};
\path (7.0,-0.7) node (v21) [label=below: $v_{2,1}$] {};
\path (8.0,0.0) node (x3) [label=below: $x_3$] {};
\path (8.0,0.85) node (v31) [label=above: $v_{3,1}$] {};
\path (9.0,0.0) node (s) [label=below: $s$] {};
\draw (r)
  -- (x1)
  -- (x2)
  -- (x3)
  -- (s);
\draw (r)
  -- (v11)
  -- (v12)
  -- (x2);
\draw (x1)
  -- (v21)
  -- (x3);
\draw (x2)
  -- (v31)
  -- (s);
\draw (x1) -- (v11);
\draw (x1) -- (v12);
\draw (x2) -- (v21);
\draw (x3) -- (v31);

\end{tikzpicture}
\end{center}

\caption{An overlapping fan graph (left) of diameter 4; fans $F_1$
and $F_3$ are same-sided (upper) fans, while $F_2$ is a lower fan,
opposite-sided from $F_1$ and $F_3$. An overlapping fan graph with
unpleasant (center: $v_{1,3} = v_{2,1} = v_{3,1}$) and pleasant
(right:  relabelled) shortest $rs$-paths.} \label{3fangraphs}
\end{figure}
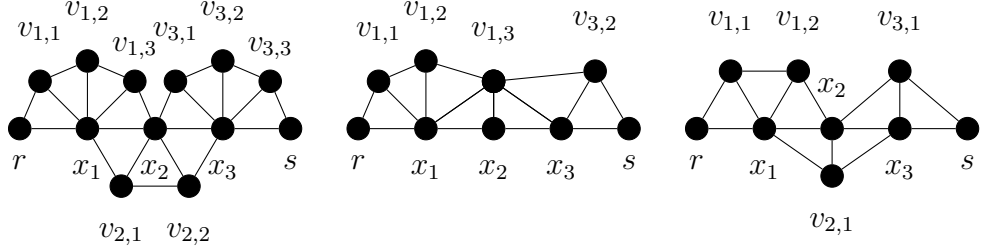

The graph $G$ is an {\it overlapping fan graph} if  the following
three conditions are satisfied: 
\begin{itemize}
    \item  for every $1\leq i \leq d-1$, there is a subgraph $F_i$
    which is  an $x_{i-1}x_{i+1}$-fan centered on $x_i$,
    \item for every $1\leq i \leq d-2$, $F_i$ and $F_{i+1}$ are
    either opposite-sided or same-sided, and
    \item $G$ is the union of the subgraphs $F_i$ for $1\leq i \leq
    d-1$.
\end{itemize}
If we agree in calling $F_1$ an {\it upper} fan, then all further
fans of an overlapping fan graph can be classified into upper or
{\it lower} (opposite-sided from upper)
--- see Figure \ref{3fangraphs}.

Notice that, in general, the description of a graph as an
overlapping fan graph, may be done using different paths $P$
(see the examples in the center and right of Figure \ref{3fangraphs}). The
path $P$ used to describe $G$ as an overlapping fan graph is
called the \textit{spine} of $G$.

In an overlapping fan graph,
$|F_i^\pr\cap F_{i+3}^\pr|=0$; while $|F_{i-1}^\pr\cap
F_{i+1}^\pr|\le 1$, with equality if and only if $k_i=1$. Notice
that  we can always choose the spine $P$ so that $|F_{i-1}^\pr\cap
F_{i+1}^\pr|=0$ by swapping the names of vertices $x_i$ and
$v_{i,1}$, changing the fans $F_{i-1}$, $F_i$, and $F_{i+1}$ from
being same-sided to $F_i$ being opposite-sided from $F_{i-1}$ and
$F_{i+1}$.  Such a choice of path $P$ is called {\it pleasant}
(see Figure \ref{3fangraphs}).

For an internal vertex $x_i$ of the spine of an overlapping fan graph $G$, we let
$A_{x_i}$ be the set of vertices of $F_i^\pr$ that are in no other fan of $G$.  If $A_{x_i}=\emptyset$ then
$k_i=1$ and $v_{i,1}\in F_{i-1}^\pr$ or $F_{i+1}^\pr$;  or $k_i=2$ and $v_{i,1}\in F_{i-1}^\pr$ and    $v_{i,2}\in F_{i+1}^\pr$. In the former let $e_{x_i}$ be the edge $x_{i-1}v_{i,1}$ or $v_{i,1}x_{i+1}$ respectively, and in the latter let   $e_{x_i}=\{v_{i,1},v_{i,2}\}$. The following fact will be used in Section \ref{otherroots}.
\begin{clm}\label{clm:e}
If $A_{x_i}$ is  empty (non empty) then  $G-e_{x_i}$  ($G-A_{x_i}$) is the union of two overlapping fan graphs each one  with $x_i$ as simplicial vertex and no other vertex in common.
\end{clm}

A 2-path of diameter 1 is just a path on two vertices.
In this case, its {\it spine} is the graph itself.
For larger diameter we have the following lemma.

\begin{lem}\label{2PathStruct}
A graph $G$ of $\diam(G)\ge 2$ is a 2-path if and only if it is an overlapping fan graph.
\end{lem}

\Pf
An overlapping fan graph is certainly a 2-path.

Let $G$ be a 2-path with simplicial vertices $r$ and $s$ and diameter at least 2.
The 2-path on 4 vertices is a fan, and hence an overlapping fan graph, 
so we assume that $G$ has at least 5 vertices.
Let $G^\pr=G-s$, with simplicial vertices $r$ and $s^\pr$.
Since $G^\pr$ is a 2-path, by induction it is also an overlapping fan graph.

If $\diam(G)>\diam(G^\pr)$ then the inclusion of $s$ creates a new fan centered on $s^\pr$.
Otherwise, the inclusion of $s$ extends the last fan of $G^\pr$.
In both cases, then, $G$ is an overlapping fan graph.
\pf

Recall that if $S$ is a set of vertices of $G$ then $G-S$ denote
the subgraph of $G$ induced by $V(G)-S$. In an analogous way, if
$F$ is a subgraph, we let $G-F$ denote the subgraph of $G$ induced
by $V(G)-V(F)$.

With respect to pebbling configurations, we define an {\it empty
vertex} (or {\it zero}) to be a vertex with no pebbles on it.  A
{\it big vertex} has at least two pebbles on it; of course, in an
$r$-unsolvable configuration, every path from a big vertex to the
root $r$ must contain at least one zero.  A {\it huge vertex} $v$
has at least $2^{\dist(v,r)}$ pebbles on it; of course, no
$r$-unsolvable configuration has a huge vertex.  The {\it cost} of
a pebbling solution $\s$ is the number of pebbles lost during the
pebbling steps of $\s$, plus one for the pebble that reaches $r$
--- we denote this by $\cost(\s)$. A {\it cheap} $r$-solution is an
$r$-solution of cost at most $2^{\ecc(r)}$, where $\ecc(r)=\ecc_G(r)$ 
is the eccentricity of $r$ in $G$.

The $t$-{\it pebbling number} $\pi_t(G)$ is the minimum number $s$
so that every configuration of size $s$ is $t$-fold solvable
(i.e., can place $t$ pebbles on any root).  The $t$-pebbling
number is related to the fractional pebbling number, which
measures the limiting average cost of repeated solutions; i.e.
$\lim_{t\rightarrow\infty}\pi_t(G)/t$.  It is also used as a
powerful inductive tool for computing the pebbling number.  The
following theorem was proven in \cite{HeHeHu}.

\begin{thm}\label{Diam2}
\cite{HeHeHu}
If $G$ is a graph of diameter 2 then $\pi_t(G)\le \pi(G)+4t-4$.
\end{thm}

In what follows we outline the key lemmas and ideas of our proof of
the pebbling number for semi-2-trees.
In Section \ref{CheapLem} we introduce the Cheap Lemma, a powerful
mechanism used in tandem with $t$-pebbling techniques.
Section \ref{2Paths} is devoted to 2-paths, which form the base step
of our induction argument for semi-2-trees in Section \ref{Semi2Trees}.
We finish with various remarks for further progress in Section \ref{Remarks}.

%
%
\section{The Cheap Lemma}\label{CheapLem}

We begin by introducing the Cheap Lemma, which we believe is a 
useful tool of independent interest.  First we develop a general
framework for some key ideas.

Fix a root $r$ in a graph $G$. 
We say that a pebbling step from $u$ to $v$ is {\it greedy} if 
$\dist(v,r)<\dist(u,r)$.
Furthermore, an $r$-solution $\s$ is {\it greedy} if each of its pebbling 
steps is greedy, and a configuration $C$ is {\it greedy} if it has a
greedy $r$-solution.
Finally, $G$ is {\it greedy} if every configuration of size at least
$\pi(G,r)$ is greedy.
(If $r$ needs to be specified, we'll use the term $r$-{\it greedy}.)

Given $\s$, let $G_\s$ denote the subgraph of edges of $G$ that are
traversed by the pebbling steps of $\s$, oriented by the direction 
of travel (bi-directed edges are allowed).
We say that $G_\s$ is {\it acyclic} if it contains no directed cycle.
The $r$-solution $\s$ is called {\it minimal} if no subset of its pebbling steps solves $r$; it is {\it minimum} if no $r$-solution uses fewer steps.
A well-known lemma of great use is the No-Cycle Lemma of \cite{BCCMW}.

\begin{lem}\label{NCL}
{\bf (No-Cycle Lemma)}
If $\s$ is a minimal $r$-solution of a configuration on $G$
then $G_\s$ is acyclic.
\end{lem}

Because of the No-Cycle Lemma, we see that every tree is greedy.
In particular, if $T$ is a breadth-first-search spanning tree of $G$, 
rooted at $r$, then $T$ is an example of an $r$-greedy spanning
subgraph of $G$ preserving distances to $r$.
Hence any configuration of size at least $\pi(T,r)$ on $G$ has
a greedy solution.
Indeed, more can be said.
Our main point will be that minimal greedy solutions are cheap,
which we will show by using weight functions.  
We say that a configuration is {\it cheap} if it has a cheap solution.

\begin{lem}\label{Cheap}
{\bf (Cheap Lemma)}
Given a graph $G$ with root $r$, let $G^*$ be an $r$-greedy spanning
subgraph of $G$ preserving distances to $r$.
Then any configuration on $G$ of size at least $\pi(G^*,r)$ is cheap.
\end{lem}

\Pf
For a vertex $v$ define the weight function $w(v)=2^{-\dist(v,r)}$;
let the {\it weight} of a configuration $C$ be $w(C)=\sum_v C(v)w(v)$.
Note that the configuration with a single pebble on $r$ has weight 1.

Suppose that $C$ is a configuration on $G$ of size at least $\pi(G^*,r)$.
Let $\s$ be a minimal greedy $r$-solution from $C$.
Denote by $C_\s$ the configuration on $G^*$ using only the pebbles of $C$
that are used by $\s$.
Then $\cost(\s)=|C_\s|$.

For any configuration $C^\pr$, let $C^{\ppr}$ be a configuration that
results from making one greedy pebbling step.  Then $w(C^\pr)=w(C^\ppr)$.
Applied iteratively to $C_\s$, this means that $w(C_\s)=1$.

Now $w(C_\s)=\sum_v C_\s(v)w(v)\ge \sum_v C_\s(v)2^{-\ecc(r)}$,
and so $\cost(\s) = |C_\s| = \sum_v C_\s(v) \le 2^{\ecc(r)}w(C_\s)
= 2^{\ecc(r)}$.
\pf

The pebbling number for a rooted tree $(T,r)$ was first derived in 
\cite{Chun}, using the notion of its {\it maximum $r$-path partition} $\cP$.
One can compute such a thing iteratively as follows.
Beginning with $F=T$, $W=\{r\}$, and $\cP=\mt$, we choose a
longest path $P$ in $F$ having one endpoint in $W$.
Then we add $P$ to $\cP$, add its vertices to $W$, remove
its edges from $F$, and repeat.

\begin{thm}\label{trees}\cite{Chun}
Let $\cP=\{P_1,\ldots,P_k\}$ be a maximum $r$-path partition of a 
rooted tree $(T,r)$, with each $P_i$ having length (number of edges) $a_i$.
(By construction, $a_i\ge a_{i+1}$ for $1\le i<k$.)
Then $\pi_t(T,r) = (t2^{a_1}-1)+\sum_{i=2}^k (2^{a_i}-1)+1 
= t2^{a_1}+\sum_{i=2}^k 2^{a_i}-k+1$.
\end{thm}

The pebbling number $\pi_t(T)$ is given by choosing $r$ to be a leaf
of a longest path of $T$.
We say that a configuration $C$ is $t$-{\it extremal} for a rooted tree 
$(T,r)$ if the following holds.
Let $\cP=\{P_1,\ldots,P_k\}$ be a maximum $r$-path partition of $(T,r)$
with each $P_i$ having leaf endpoint $v_i$.
Then $C(v_1)=t2^{a_1}-1$, $C(v_i)=2^{a_i}-1$ for $2\le i\le k$, and
$C(v)=0$ otherwise.
The proof of the lower bound in Theorem \ref{trees} involves showing
(by induction) that such a configuration is not $t$-fold $r$-solvable.

For a 2-path $G$ with simplicial root $r$, we denote by $T^*(G,r)$ any
spanning tree of $G$, rooted at $r$, that includes the spine of $G$ and
all fan vertices as leaves, each one adjacent to its neighbor in the spine closest to $r$. 
Notice that $T^*(G,r)$ is an $r$-greedy spanning subgraph of $G$ preserving distances to $r$.

For a 2-path $G$ with  simplicial
vertex $r$,  root eccentricity $d$, and  with $n$ vertices, we define the functions $p_t(G,r)=t2^d+n-2d$ (suppressing
$t$ when $t=1$) and $q(G,r)=2^d+n-d-1$.
Note that $p(G,r)<q(G,r)<p_2(G,r)$ when $1<d$.

\begin{cor}\label{tCheap}
Let $G$ be a 2-path with simplicial vertex $r$ and diameter $d$.
If $C$ is a configuration of size
at least $q(G,r)+(t-1)2^d$ then $C$ has $t$ distinct cheap $r$-solutions.
\end{cor}

\Pf
For $t=1$ this follows from the Cheap Lemma \ref{Cheap} and Theorem 
\ref{trees} because for $T^*=T^*(G,r)$ we have $\pi(T^*,r)=q(G,r)$.
The general statement follows by induction on $t$.
\pf


The following two lemmas about pebbling in trees will be used in Section \ref{otherroots}.

\begin{lem}\label{TreeLemma}
Let $T$ be a tree with diameter $d=\diam(T)$, $r^*$ and $r$ be vertices with $\ecc(r)<\ecc(r^*)=d$.
Let $P^*$ be a path $v_0v_1\cdots v_d$ with $v_0=r^*$ and $v_d=s^*$, labeled so that $\dist(r,s^*)\le\dist(r,r^*)=\ecc(r)$.
Denote by $P$ the path from $r$ to $r^*$, and set $P^*\cap P=v_0\cdots v_{h^\pr}$.
Define $\overline{h}=d-h^\pr$.
Then $\pi_t(T,r) \le \pi_t(T,r^*) - t(2^d-2^{\ecc(r)}) + 2^{\overline{h}} - 1 \le \pi_t(T,r^*) - 2^{d-2}$.
\end{lem}

\Pf
Let $\cP^*$ be a maximum path partition of $T$ with root $r^*$.
Define $P^*_0=P^*$, $P^*_1$, $\ldots$, $P^*_k$ to be the sequence of paths of $\cP^*$ that are used sequentially while traveling from $r^*$ to $r$ in $P$, and set $d^*_i=\len(P^*_i)$ for each $0\le i\le k$ (so $d^*_0=d$).
Next define $P_i^\pr=P\cap P^*_i$, with $h^\pr_i=\len(P_i^\pr)$ and $\overline{h}_i=d^*_i-h^\pr_i$ (so $h^\pr_0=h^\pr$ and $\overline{h}_0=\overline{h}$).
Notice that $\ecc(r)=\sum_{i=0}^kh^\pr_i$ and $\overline{h}\le d/2$.

Denote by $\cP$ the maximum path partition of $T$ with $r$ as root.
We will use the following facts in the calculations below.
\begin{itemize}
\item 
The longest path in $\cP$ is $P$.
\item
In the component of the tree $T-P$ that contains the path $\hat{P_i}=P^*_i-P_i^\pr$, the longest path is $\hat{P_i}$.
\end{itemize}
From these it follows that each $\hat{P_i}\in\cP$ and, subsequently, that $\cP^*-\{P^*_0,\ldots,P^*_k\}=\cP-\{P,\hat{P_0},
\ldots,\hat{P_k}\}$.
Now, by converting $\cP^*$ to $\cP$, we find that
\begin{align*}
\pi_t(T,r)
& = \pi_t(T,r^*) - \left[(t2^d-1) + \sum_{i=1}^k(2^{d^*_i}-1)\right] + \left[(t2^{\ecc(r)}-1) + \sum_{i=0}^k(2^{\overline{h}_i}-1)\right]\\
& \le \pi_t(T,r^*) + t2^{\ecc(r)} - \left(t2^d-2^{\overline{h}_0} \right) - \left[\sum_{i=1}^k \left(2^{d^*_i}-2^{\overline{h}_i} \right) \right] - 1\\
& \le \pi_t(T,r^*) - t(2^d-2^{\ecc(r)}) + 2^{\overline{h}} - 1\\
& \le \pi_t(T,r^*) - t2^{d-1}+2^{\lfloor d/2\rfloor} - 1\\
& \le \pi_t(T,r^*) - 2^{d-2}\ .
\end{align*}
\pf

\begin{lem}\label{l:tree3} 
Let $e=xy$ be a non pendant edge of a tree $T$ and assume that $\ecc(x)\geq \ecc(y)$. 
If $T'$ is the tree obtained by subdividing the edge $e$ with a new vertex $r$, then $\pi_t(T',r)= \pi_t(T,x)+2^{a}$, where $a$ is the eccentricity of $x$ in the connected component of $T-y$ that contains $x$ (thus, $a+\ecc(x)\leq \diam(T)$).
\end{lem}

\Pf 
Define $x^\pr$ to be the vertex having $\dist_T(x,x^\pr)=\ecc_T(x)$, and denote the $xx^\pr$-path by $P$.
Because $\ecc_T(y)\le \ecc_T(x)$ we know that $y\in P$.
Let $x^\ppr$ be a vertex having $\dist_{T-y}(x,x^\ppr)=a$, with $xx^\ppr$-path $Q$, and note that $a<\ecc_T(x)$.

Now observe that $\dist_{T^\pr}(r,x^\pr)=\dist_T(x,x^\pr)$ (witnessed by the $rx^\pr$ path $P^\pr$) and $\dist_{T^\pr}(r,x^\ppr)=\dist_T(x,x^\ppr)+1$ (witnessed by the $rx^\ppr$ path $Q^\pr$).
This means that the only changes from the maximum path partition of $T$ with root $x$ to the maximum path partition of $T^\pr$ with root $r$ are that the longest path $P$ from $x$ in $T$ becomes the longest path $P^\pr$ from $r$ in $T^\pr$, and the longest path $Q$ from $x$ in $T-y$ becomes the longest path $Q^\pr$ from $r$ in $T^\pr-y$.
Hence we have $\pi_t(T,x)=t2^{\ecc_T(x)}+2^a+F(x)$, for some $F(x)$, and $\pi_t(T^\pr,r)=t2^{\ecc_{T^\pr}(r)}+2^{a+1}+F(x) = \pi_t(T,x)+2^{a+1}-2^a$.
\pf

%
%
\section{ 2-Paths}\label{2Paths}

In this section we calculate  a $\pi_t(G,r)$ for $r$ a simplicial vertex of a 2-path $G$.

\subsection{The Lower Bound}\label{lowerbound}

We now present some general removal techniques for finding lower bounds 
that may also be of independent interest.
For a vertex $v$, define its {\it open neighborhood} $N(v)$ to be the set of
vertices adjacent to $v$, and its {\it closed neighborhood} $N[v]=N(v)\cup\{v\}$.
Also, for a set of vertices $A$ write $N(A)=\cup_{v\in A}N(v)$.
Along the lines of the definition of twin vertices,
for a non-root vertex $y$ we say that $y$ is a {\it junior sibling of} $x$ 
(or, more simply, {\it junior to} $x$) if $N(y)\sse N[x]$,
and that $y$ is a {\it junior} if it is junior to some vertex $x$.

\begin{lem}{\bf (Junior Removal Lemma)}\label{JuniorRemoval}
Given the rooted graph $(G,r)$ with configuration $C$, 
suppose that $y$ is a junior with $C(y)=0$.
Then $C$ is $t$-fold $r$-solvable if and only if $C$ restricted to $G-y$ is $t$-fold 
$r$-solvable in $G-y$.
\end{lem}

\Pf
Sufficiency is obvious, so we only prove necessity.
Suppose that $\s$ is an $r$-solution from $C$ that uses $y$.
Let $y$ be junior to some vertex $x$.
Construct $\s^\pr$ from $\s$ by replacing every pebbling step $(u,y)$
with $(u,x)$ and every pebbling step $(y,v)$ with $(x,v)$.
Then $\s^\pr$ $t$-fold solves $r$ as well.
\pf

We say that a set of vertices $W$ is a {\it wart} if it is a component
of $G-X$ for some clique cutset $X$, where by {\it clique} we mean complete subgraph.

\begin{lem}{\bf (Wart Removal Lemma)}\label{WartRemoval}
Given the rooted graph $(G,r)$ with configuration $C$, 
suppose that $W$ is a wart of $G$ not containing $r$ and that
$C(w)\le 1$ for every $w\in W$.  
Then $C$ is $t$-fold $r$-solvable if and only if $C$ restricted to $G-W$ is 
$t$-fold $r$-solvable in $G-W$.
\end{lem}

\Pf
Sufficiency is obvious, so we only prove necessity.
We show that no minimum $r$-solution from $C$ uses $W$.

Suppose instead that $\s$ is a minimum $r$-solution that uses $W$.
Let $X$ be a clique cutset that witnesses the wart $W$, and let $u$ be a vertex of $X$ having a pebbling step into $W$.
Because $\s$ is minimum, there is a vertex $v\in X$ that receives
a pebble from $W$ and that is different from $u$.
By replacing those two pebbling steps by the single step from $u$ to $v$
we find an $r$-solution with fewer steps, a contradiction.
\pf

Let $G$ be a 2-path with simplicial root $r$, pleasant path $P$, and configuration $C$.  
For a given $t$ we say that $C$ is $t$-{\it extremal for} $r$
(simply, {\it extremal} if $t=1$) if there is a $I$-saturating matching $M$ from
the internal spine vertices $I=\{x_1, \ldots, x_{d-1}\}$ to the fan
vertices $\{v_{i,j}\}$ such that $C(x_d)=t2^d-1$, $C(r)=0$, $C(M)=0$,
and $C(v)=1$ otherwise.  Notice that $|C|=p_t(G,r)-1$.

If a configuration $C$ on $G$ is $t$-fold $r$-solvable if and only if
$C_H$ is $t$-fold $r$-solvable on the subgraph $H\subset G$,
then we say that $G$ $t$-{\it fold} $r$-{\it reduces to} $H$ {\it for} $C$.
If $C$, $t$ and $r$ are clear from the context we just write {\it reduces}.

\begin{lem}{\bf (Extremal Lemma)}\label{extremal}
If $C$ is $t$-extremal for the simplicial root $r$ of a 2-path $G$
then $C$ is not $t$-fold $r$-solvable.
Moreover, by using Lemmas \ref{JuniorRemoval} and \ref{WartRemoval}
(repeatedly removing juniors and warts) 
$G$ reduces to its spine, the path $P_d$, where $d=diam(G)$.
\end{lem}

\Pf
We use induction on $d$.
The result is trivial for $d=1$.
For $d>1$ we suppose that $C$ is $t$-fold $r$-solvable
and let $\s$ be a $t$-fold $r$-solution.
Write $y_i=v_{i,j_i}$ for the neighbor of $x_i$ in $M$ and
let $\ell$ be the smallest index $i$ such that $y_i$ is a junior.
This exists because if $y_i$ is not a junior then either $y\in F_{i-1}\cap F_i$
or $y\in F_i\cap F_{i+1}$ (it is a {\it fan intersection}), and there are
more fans than fan intersections.
Set $y=y_{\ell}$ and $x=x_{\ell}$.  
Then $y$ is junior to $x$ and so, by Lemma \ref{JuniorRemoval}, 
$C$ is $t$-fold $r$-solvable in $G-y$.

Furthermore, let $j^+$ be the maximum $j$ such that $v_{\ell,j}\in F_{\ell}-F_{\ell +1}$.
If $j_{\ell}+1\le j^+$ then $\{v_{\ell,j_{\ell}+1}\}$ is a wart in $G-y$, and so Lemma
\ref{WartRemoval} says that we can remove it.
Once we do, $\{v_{\ell,j_{\ell}+2}\}$ becomes a wart, and so on, until
all the vertices $v_{\ell,j}$ with $j_{\ell}<j\le j^+$ have been removed.
Then the graph $G_{\ell +1}=\cup_{i>\ell}F_i$ is a 2-path, with the restriction,
$C_{\ell +1}$, of $C$ to $G_{\ell +1}$ being $2^{\ell} t$-extremal for $x_{\ell}$.
By induction, $C_{\ell +1}$ is not $2^{\ell} t$-fold $x_{\ell}$-solvable and
$G_{\ell +1}$ can be reduced to the path $P_{d- \ell}$.

Similarly, let $j^-$ be the minimum $j$ such that $v_{\ell,j}\in F_{\ell}-F_{\ell -1}$.
If $j_{\ell}-1\ge j^-$ then the warts $\{v_{\ell,j}\}$ for $j^-\le j<j_{\ell}$ can be successively
removed, leaving the 2-path $G^{\ell}=\cup_{i\le \ell}F_i$.
Since $C$ is $t$-fold $r$-solvable and $x_{\ell}$ is a cut-vertex of $G-y$,
all the pebbles of $G_{\ell +1}$ used by $\s$ must pass through $x_{\ell}$.
But because $C_{\ell +1}$ is not $2^{\ell} t$-fold $x_{\ell}$-solvable, the most number
of pebbles that can reach $x_{\ell}$ is $2^{\ell} t-1$.
After placing as many pebbles as possible on $x_{\ell}$ from $G_{\ell+1}$, the
resulting configuration $C^{\ell}$ is a subconfiguration of a configuration
$\hat{C}^{\ell}$ that is $t$-extremal for $r$ on $G^{\ell}$.
By induction, $\hat{C}^{\ell}$ is not $t$-fold $r$-solvable, a contradiction.
Also, $G^{\ell}$ can be reduced to the path $P_{\ell}$, which reduces $G$ to the
path $P_d$.
\pf

\begin{cor}\label{DdLB}
If $r$ is a simplicial vertex of a 2-path  $G$ then $\pi_t(G,r)\ge p_t(G,r)$.
\pf
\end{cor}

\subsection{The Upper Bound}\label{upperbound}

We first note that a diameter two 2-path $G$ is Class 0.
Indeed, the following lemma is a corollary of the Class 0 characterization for diameter two
graphs from \cite{ClHoHu} that shows that $\pi(G)=n$ in this case and 
the $t$-pebbling bound of \cite{HeHeHu} that states $\pi_t(G)\le 
\pi(G)+4t-4$ for all diameter two graphs.
Equality comes from Corollary \ref{DdLB}.
The diameter one case is from \cite{HeHeHu} also.

\begin{lem}\cite{HeHeHu}\label{D2Class0}
If $G$ is a 2-path on $n$ vertices with diameter $d \leq 2$ then
$\pi_t(G)=t2^d+n-2d$. 
\end{lem}

\begin{thm}\label{PathSimplicial}
Let $G$ be a 2-path on $n$ vertices with simplicial root vertex $r$
having eccentricity $d$, and configuration $C$.
If $|C|\ge p(G,r)$ then $C$ is $r$-solvable.
\end{thm}

\Pf 
When $d\leq 2$, the result is taken care of by Lemma
\ref{D2Class0}. So we will assume that $d>2$ and use induction.
Suppose that $|C|=p(G,r)$ and let $P=r,x_1,\ldots,x_{d-1},s$ be
a pleasant shortest $rs$-path between the two simplicial vertices
of $G$. Write $x_0=r$ and $x_d=s$ and label $G$ by its fan graph labeling, so that
$V(F_i)=\{x_{i-1},x_i,x_{i+1}, v_{i,1},\ldots,v_{i,k_i}\}$ and
$Q_i$ is the path $x_{i-1} ,v_{i,1},\ldots,v_{i,k_i},x_{i+1}$. Let
$\Gp$ be the restriction of $G$ to the $\np$ vertices of
$\cup_{i\ge 2}V(F_i)$, with $\Cp$ denoting the restriction of $C$
to $\Gp$. We further use the abbreviations $C_1=C(F_1)$ and
$n_1=|V(F_1)|$. Notice that $\diam(\Gp)=d-1$, so that the Theorem
holds for $\Gp$.
Define $\f=1$ ($0$) if $F_2$ is same-sided (opposite-sided) as $F_1$.

If $C(x_1)\ge 1$, $C(x_2)\ge 2$, or $C(v_{1,j})\ge 2$ for some $j$ 
(either [$\f=1$ and $j=k_1$] or not), then we can place a pebble on $x_1$. 
If $|\Cp|-(1,2,2,0)\ge p(\Gp,x_1)$, where the coordinates correspond, 
in order, to the four cases above (first two cases plus two sub-cases
of the third case), then we can place another pebble on $x_1$,
and then one on $r$.
Otherwise, $|\Cp|-(1,2,2,0)\le p(\Gp,x_1)-1$. 
That is, $|\Cp|\le [2^{d-1}+\np-2(d-1)]+(0,1,1,-1)$. 
Thus $|C_1|\ge |C|-|\Cp|+(1,2,2,0)\ge 2^{d-1}+(n_1-2-\f)-2+(1,1,1,1) =
n_1+(2^{d-1}-3-\f)\ge n_1$, which means by Lemma \ref{D2Class0} that we
can solve $r$.

On the other hand, if $C(x_1)=0$, $C(x_2)\le 1$, and $C(v_{1,j})\le 1$
for all $j$, then $C(\{r, v_{1,1}, \ldots,$ $v_{1,k_1-1}\})\le k_1-1$.
Here we define $\th$ to be the number of zeros in 
$\{v_{1,1},\ldots,v_{1,k_1},x_2\}$, so that $|C_1|=n_1-2-\th$, and set $\thp$ 
to be the number of those zeros other than $x_2$ (i.e. $\th-\thp=1-C(x_2)$).
Now we have
\begin{align*}
|\Cp|
 & \ge |C|-|C_1|+C(x_2)\\
 & = (2^d+n-2d)-(n_1-2-\th)+C(x_2)\\
 & = (2)2^{d-1} + (\np-2-\f)-2d+2+\th+C(x_2)\\
 & = [(2)2^\dpr+\np-2\dpr] + [C(x_2)+\th-2-\f]\\
 & = p_2(\Gp,x_1) + [\thp-1-\f]\ .\\
\end{align*}
If $\thp-1-\f\ge 0$ then $|\Cp|\ge p_2(\Gp,x_1)>q(\Gp,x_1)$, which means, by
Corollary \ref{tCheap}, that we can place one pebble on $x_1$ cheaply.
Because the remaining configuration (after solving $x_1$ cheaply)
has at least $p_2(\Gp,x_1)-2^\dpr = p(\Gp,x_1)$ pebbles, 
induction places a second pebble on $x_1$.	
Then we move one to $r$.

Otherwise, we have $\thp-1-\f<0$, which means that $\thp\le\f$.  If
$\thp=0$, that is $C(v_{1,j})=1$ for all $j$, then we will show
that it is possible to place two pebbles on $x_2$, from which we
solve $r$ by moving pebbles from $x_2$ along $Q_1$. Indeed, this
is so if $C(x_2)=1$ and $Q_2$ has a big vertex, or if $Q_2$
contains either a vertex with four pebbles or two big vertices, so
we assume otherwise. In this case, we have $|C((F_1\cup
F_2)-\Gpp)|\le |V((F_1\cup F_2)-\Gpp)|$, where $\Gpp$ is the restriction
of $G$ to the $\npp$ vertices of $\cup_{i\ge 3}V(F_i)$. For the
restriction $\Cpp$ of $C$ to $\Gpp$, this implies that
\begin{align*}
|\Cpp|
 & = |C|- |C((F_1\cup F_2)-\Gpp)|\\
 & \ge 2^d + \npp - 2d\\
 & = [(2)2^{d-2} + \npp - 2(d-2)] + [2^{d-1}-4]\\
 & \ge p_2(\Gpp,x_2)\ ,
\end{align*}
since $d\ge 3$.
As before, since $p_2(\Gpp,x_2)>q(\Gpp,x_2)$ and $p_2(\Gpp,x_2)-2^\dppr
=p(\Gpp,x_2)$, we can place one pebble on $x_2$ cheaply, followed by
a second pebble on $x_2$.

We are left now with the final case (since $\thp\le\f\le 1$)
in which $\thp=1$ (exactly one $v_{1,j}$
is empty), which means that $\f=1$ ($F_1$ and $F_2$ are same-sided,
so that $v_{1,k_1}=v_{2,1}$).

If $v_{1,k_1}$ is not empty then $k_1\ge 2$, and so
\begin{align*}
|\Cp|
 & = |C|-(k_1-2)\\
 & = (2)2^{d-1} + (n-k_1) - 2(d-1)\\
 & = p_2(\Gp,x_1)\\
 & > q(\Gp,x_1)\ .
\end{align*}
As above, this means, by Corollary \ref{tCheap} and induction, that we 
can place two pebbles on $x_1$, and hence one on $r$.

If instead $v_{1,k_1}$ is empty then set $\Gh = \Gp-x_1$ and
$\Ch=C(\Gh)$, so that
\begin{align*}
|\Ch|
 & = |C|-(k_1-1)\\
 & = (2)2^{d-1} + (n-k_1-1) - 2(d-1)\\
 & = p_2(\Gh,v_{1,k_1})\ .
\end{align*}
Again, this means that we can place two pebbles on $v_{1,k_1}$,
and hence one on $r$ (via $Q_1$).

This completes the proof.
\pf

\begin{cor}\label{tPebbSimplicial}
If $r$ is a simplicial vertex of a 2-path $G$ then $\pi_t(G,r)=p_t(G,r)$.
\end{cor}

\Pf 
The lower bound was stated in Corollary \ref{DdLB}.
The upper bound for $t=1$ follows from Theorem \ref{PathSimplicial}.
If $t>1$, then for any configuration $C$ of size $p_t(G,r) =
p_2(G,r)+(t-2)2^d > q(G,r)+(t-2)2^d$, we can place $t-1$ pebbles 
on $r$, each cheaply, by Corollary \ref{tCheap}.
The remaining configuration has at least $p_t(G,r)-(t-1)2^d =
p(G,r)$ pebbles, from which we can place the $t^{\rm th}$ pebble
on $r$ by Theorem \ref{PathSimplicial}.
\pf

%
%
\section{Pebbling number of Semi-2-Trees}\label{Semi2Trees}

We define the {\it skeleton} $T$ of a semi-2-tree $G$ to be the union
of the spines of its blocks; it is a geodesic tree spanning all of the
simplicial vertices of $G$.
Let $e(T)$ denote the number of edges of $T$, $b(G)$ denote the
number of blocks of $G$, and for a simplicial vertex or cut-vertex $r$ 
and positive integer $t$ define $p_t(G,r)=\pi_t(T,r)+(n-1)+b(G)-2e(T)$
(suppressing $t$ when $t=1$).
Notice that this matches the corresponding formula for 2-paths because $b=1$ and $T$ is a path.
In addition, we have $p_t(G,r)=p_{t-1}(G,r)+2^{\ecc_G(r)}$ because of Theorem \ref{trees}.
We also define $q(G,r)=\pi(T,r)+n-e(T)-1$; note that $q(G,r)=\pi(T^*,r)$, where $T^*$ is a spanning tree of $G$, rooted at $r$, that contains its skeleton and all its fan vertices as leaves, each one adjacent to its neighbor in the skeleton closest to $r$.
Notice that $T^*$ is an $r$-greedy spanning tree of $G$ preserving distances to $r$.

\subsection{Simplicial or Cut-vertex Roots}\label{SimplicialOrCut}

We begin with another consequence of the Cheap Lemma, generalizing
Corollary \ref{tCheap}.  The proof is similar and is left to the reader.

\begin{cor}\label{tCheapTree}
Let $r$ be a simplicial vertex or cut-vertex with eccentricity $d$ of a semi-2-tree $G$. If $C$ is a configuration 
of size at least $q(G,r)+(t-1)2^d$ then $C$ has $t$ distinct cheap $r$-solutions.
\pf
\end{cor}

For a tree $T$ with maximum $r$-path partition $\cP=\{P_1,\ldots,P_k\}$, 
each $P_i$ having length $a_i$ (sorted so that $a_i\ge a_{i+1}$),
let $C_T$ be its $t$-extremal configuration for $r$.

For a semi-2-tree $G$, call a vertex of the skeleton $T$  {\it internal} if it is
not a simplicial vertex or cut-vertex, and let $M$ be any $I$-saturating matching
from the internal vertices $I$ to the fan vertices of $G$. For a simplicial or cut vertex $r$ of $G$, define the configuration $C$ by $C(T)=C_T$, $C(M)=0$, and $C(v)=1$
otherwise --- such a configuration we call $t$-{\it extremal for} $r$.
Note that $|C|=p_t(G,r)-1$.

As in the proof of the Extremal Lemma \ref{extremal}, we can use the Removal 
Lemmas \ref{JuniorRemoval} and \ref{WartRemoval} to prove that $G$ reduces to $T$ for $C$ and obtain the following more general extremal lemma, which we leave to the reader.

\begin{lem}\label{extremaltree}
If $C$ is $t$-extremal for the simplicial or cut-vertex root $r$ of a semi-2-tree $G$
then $C$ is not $t$-fold $r$-solvable.
Moreover, by using Lemmas \ref{JuniorRemoval} and \ref{WartRemoval},
$G$ can be reduced to its skeleton $T$.
\pf
\end{lem}

Now we state and prove the solvability theorem in this case.

\begin{thm}\label{TreeSimplicial}
Let $G$ be a semi-2-tree on $n$ vertices with simplicial or cut-vertex
$r$ and configuration $C$.
If $|C|\ge p(G,r)$ then $C$ is $r$-solvable.
\end{thm}

\Pf

We use induction on $n$ with base case $\ecc(r)=1$, which is handled
by Theorem \ref{Diam2}.
So we assume that $\ecc(r)>1$.
We may also assume that $C(r)=0$.
We consider two cases.

\begin{enumerate}
\item
$r$ is a cut-vertex.

Let $H_1,\ldots,H_k$ be the components of $G-r$, with $G_i$ induced
by $V(H_i)\cup\{r\}$; then each $G_i$ is a semi-2-tree, so that
the theorem holds for them by induction.
Let $C_i$ and $T_i$ be the restrictions of $C$ and $T$ to $G_i$,
with $n_i$ and $b_i$ counting the number of vertices and blocks of $G_i$.
If some $|C_i|\ge p(G_i,r)$ then $C_i$ solves $r$, so we assume not.
Then
\begin{align*}
|C|
 & = \sum_i |C_i|\\
 & \le \sum_i [p(G_i,r)-1]\\
 & = \sum_i [\pi(T_i,r) + (n_i-1) + b_i - 2e(T_i) - 1]\\
 & = \pi(T,r) + (n-1) + b(G) - 2e(T) - k\\
 & < p(G,r)\ ,
\end{align*}
a contradiction.
Hence some $C_i$ solves $r$.
\item
$r$ is a simplicial vertex.

Let $H$ be the block of $G$ containing $r$, with $\rp$ the other
simplicial vertex of $H$.
If $|C(H)|\ge p(H,r)$ then we solve $r$ directly on $H$.
Otherwise, we assume that $|C(H)|=p(H,r)-s$ for some $s>0$.
Recall that $p(H,r)=2^\dH+\nH-2\dH$, where $\nH=|H|$ and $\dH=\ecc_H(r)$.
Let $\Gp$ be the subgraph of $G$ induced by $(G-H)\cup\{\rp\}$, having
$\np=n-\nH+1$ vertices, $\bp$ blocks, and root eccentricity 
$\ecc_\Gp(\rp)=\dpr=d-\dH$, with $\Tp=T\cap\Gp$ and $d=\ecc_G(r)$.
Define the configuration $\Cp$ on $\Gp$ by $\Cp(\rp)=0$ and $\Cp(v)=C(v)$
for all other $v\in\Gp$.

Suppose that $s\le 2^\dH$.
Then
\begin{align*}
|\Cp|
 & = |C|-|C(H)|\\
 & = p(G,r)-p(H,r)+s\\
 & = [\pi(T,r)+(n-1)+b(G)-2e(T)]-[2^\dH+\nH-2\dH]+s\\
 & = [\pi_s(\Tp,\rp)+(\np-1)+\bp-2e(\Tp)]+[s-2^\dH+2^d-s2^\dpr]\\
 & = p_s(\Gp,\rp)+(2^\dpr-1)(2^\dH-s)\\
 & \ge p_s(\Gp,\rp)\ ,
\end{align*}
which means that we can place $s$ pebbles on $\rp$, so that now there
are $p(H,r)$ pebbles in $H$, enough to solve $r$.

Suppose that $s\ge 2^\dH$; i.e. $|C(H)|\le \nH-2\dH$.
Then
\begin{align*}
|\Cp|
 & = |C|-|C(H)|\\
 & \ge [\pi(T,r)+(n-1)+b(G)-2e(T)]-[\nH-2\dH]\\
 & = [\pi_{2^\dH}(\Tp,\rp)+(\np-1)+\bp-2e(\Tp)]\\
 & \ge p_{2^\dH}(\Gp,\rp)\ ,
\end{align*}
which means that we can place $2^\dH$ pebbles on $\rp$, 
enough to solve $r$ on $T$.
\end{enumerate}
\pf

\begin{cor}\label{tPebbSimplicialTree}
If $r$ is a simplicial vertex or cut-vertex of a semi-2-tree $G$ then $\pi_t(G,r)=p_t(G,r)$.
\end{cor}

\Pf
As in the proof of Corollary \ref{tPebbSimplicial}.
\pf

\begin{thm}\label{BestRoot1}
If $r$ is a simplicial vertex or cut-vertex of a semi-2-tree $G$ and $r^*$ is a simplicial vertex with $\ecc(r^*)=\diam(G)$ then $\pi_t(G,r)\le\pi_t(G,r^*)$.
\end{thm}

\Pf
Let $T$ be a skeleton of $G$.
Because the only term in $p_t(G,r)=\pi_t(T,r)+(n-1)+b(G)-2e(T)$ that depends on $r$ is $\pi_t(T,r)$, it follows that $\pi_t(G,r)$ is maximized precisely where $\pi_t(T,r)$ is maximized, which is well-known (\cite{Chun}) to be at $r^*$.

\pf

\subsection{Other Roots}\label{otherroots}

We begin with two more removal lemmas of general use.

\begin{lem}{\bf (Edge Removal Lemma)}\label{EdgeRemoval}
Let $r$ be a vertex of a connected graph $G$ and suppose $e$ is an
edge between two neighbors of $r$.
Then $\pi(G,r)=\pi(G-e,r)$.
\end{lem}

\Pf
Given any configuration on $V(G)=V(G-e)$, every minimal $r$-solution in 
one graph is a minimal solution in the other.
\pf

\begin{lem}\label{Components}
Let $r$ be a cut-vertex of a graph $G$, and denote the connected
components of $G-r$ by $H_1,\ldots,H_k$.  For each $i$ define the graph
$G_i$ induced by $H_i\cup\{r\}$.
Then $\pi(G,r)=1+\sum_i (\pi(G_i,r)-1)$.
\end{lem}

\Pf
The lower bound follows from the union of the individual maximum-sized
\mbox{$r$-unsolvable} configurations on $H_i$.
The upper bound follows from the pigeonhole principle.
\pf

\begin{lem}{\bf (Neighbor Removal Lemma)}\label{NeighborRemoval}
Let $r$ be a vertex of a connected graph $G$. 
Suppose that $A\sse N(r)$ such that $N(A)\sse N[r]$.
Let $\{H_1,\ldots,H_k\}$ be the connected components of $(G-r)-A$ and 
denote by $G_i$ the subgraph of $G$ induced by $V(H_i)\cup \{r\}$. 
Then $\pi(G,r)= 1+|A|+\sum_i (\pi(G_i,r)-1)=|A|+\pi(G-A,r)$.
\end{lem}

\Pf
We can remove the edges incident with $A$ by Lemma \ref{EdgeRemoval}.
Then each $v\in A$ is its own component of $G-r$.
The result follows from Lemma \ref{Components}.
\pf

Under the conditions of Lemma \ref{NeighborRemoval}, if each $(G_i,r)$ is a rooted semi-2-tree, then we say that a configuration $C$ on $G$ is {\it extremal for} $r$ if $C(x)=1$ for every $x\in A$ and each $C_{G_i}$ is extremal for $r$ on $G_i$.

\begin{figure}
\begin{center}
\begin{tikzpicture}[scale=1.2]
\tikzstyle{every node}=[draw,circle,fill=black,minimum size=2pt,inner sep=3pt]
\path (2.0,0.0) node (z) [label=below: $z$] {};
\path (3.0,0.0) node (y) [label=below: $y$] {};
\path (3.5,0.85) node (r) [label=above: $r$] {};
\path (4.0,0.0) node (x) [label=below: $x$] {};
\draw (z)
  -- (y)
  -- (x);
\draw (y)
  -- (r)
  -- (x);
\end{tikzpicture}
\end{center}
\caption{A non-semi-2-tree $G$ for which $G-x$ is a semi-2-tree.
The configuration $C(r,x,y,z)=(0,1,0,3)$ is extremal for $r$.} \label{SmallNbrRem}
\end{figure}
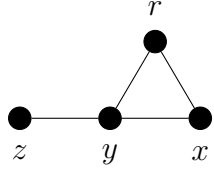

A small example of a non-semi-2-tree to which Lemma \ref{NeighborRemoval} applies is shown in Figure \ref{SmallNbrRem}.
This idea is used later in the proof of Corollary \ref{onefanold}.

A simple consequence (using the Cheap Lemma and induction) of Lemma \ref{EdgeRemoval} is the following.

\begin{cor} \label{spinepebbempty}
Let $G$ be a semi-2-tree with skeleton $T$, and suppose that $r$ is a vertex of $T$ that is not a simplicial or cut vertex of $G$. Let $A_r$ be the set of vertices of the fan centered on $r$ that are in no other fan of $G$. If $A_r$ is empty and $e_r$ is as defined on Claim \ref{clm:e},  then $\pi_t(G,r)=\pi_t(G-e_r,r)$ for all $t\ge 1$.
\pf
\end{cor}

Notice that the previous corollary allows one to calculate the pebbling number for $r$. In fact, by Claim \ref{clm:e}, $G-{e_r}$ is a semi-2-tree with $r$ a simplicial or cut vertex, then  we use %
Corollary \ref{tPebbSimplicialTree} to calculate $\pi(G,r)=\pi(G-e_r,r)$.

Anagously, a consequence (using the Cheap Lemma and induction) of Lemma \ref{NeighborRemoval} is the following.

\begin{cor} \label{spinepebb}
Let $G$ be a semi-2-tree with skeleton $T$, and suppose that $r$ is a vertex of $T$ that is not a simplicial or cut vertex of $G$.  Let $A_r$ be the set of vertices of the fan centered on $r$ that are in no other fan of $G$. If $A_r$ is non empty then $\pi_t(G,r)=\pi_t(G-A_r,r)+|A_r|$ for all $t\ge 1$.
\pf
\end{cor}

Notice that the previous corollary allows one to calculate the pebbling number for $r$. In fact, by Claim \ref{clm:e}, $G-{A_r}$ is a semi-2-tree with $r$ a simplicial or cut vertex,   then  we use
Corollary \ref{tPebbSimplicialTree} to calculate $\pi(G,r)=\pi(G-A_r,r)$.

\begin{thm}\label{BestRoot2}
Let $G$ be a semi-2-tree with skeleton $T$.
Suppose that $r$ is a vertex of $T$ that is not a simplicial or cut vertex of $G$, and let $r^*$ be a simplicial vertex of $G$ with $\ecc(r^*)=\diam(G)$.
Then $\pi_t(G,r)<\pi_t(T,r^*)$.
\end{thm}

\Pf
Let $A_r$ be the set of vertices of the fan centered on $r$ that are in no other fan of $G$. 
First assume $A_r=\emptyset$ and let $e_r$ be as in Corollary \ref{spinepebbempty}. Notice that
the skeleton of $G-e_r$ is the same $T$, while  $G-e_r$ has one block less than $G$, thus (using Corollary \ref{spinepebbempty} and Corollary \ref{tPebbSimplicialTree})
\begin{align*}
\pi_t(G,r)
& = \pi_t(G-e_r,r)\\
& = p_t(G-e_r,r)\\
& = \pi_t(T,r)+(n(G-e_r)-1)+b(G-e_r)-2e(T)\\
& = \pi_t(T,r)+(n(G)-1)+(b(G)-1)-2e(T)\\
& \leq  \pi_t(T,r^*)+(n(G)-1)+b(G)-2e(T)-1\\
& = \pi_t(G,r^*)-1\\
& < \pi_t(G,r^*)\ .
\end{align*}

Analogously, if $A_r\neq\emptyset$ then (using Corollary \ref{spinepebb})

\begin{align*}
\pi_t(G,r)
& = \pi_t(G-A_r,r)+|A_r|\\
& = p_t(G-A_r,r)+|A_r|\\
& = \pi_t(T,r)+(n(G-A_r)-1)+b(G-A_r)-2e(T)+|A_r|\\
& = \pi_t(T,r)+((n(G)-|A_r|)-1)+(b(G)-1)-2e(T)+|A_r|\\
& = \pi_t(T,r)+(n(G)-1)+b(G)-2e(T)-1\\
& \le \pi_t(G,r^*)-1\\
& < \pi_t(G,r^*)\ .
\end{align*}

\pf

Another consequence (again using the Cheap Lemma and induction) of Lemma \ref{EdgeRemoval} is the following. We say that a vertex   $r$ is not in any skeleton of $G$ when for every skeleton $T$ of $G$, $r$ is a fan vertex, i.e. $r\in V(G-T)$.

\begin{cor}\label{twofans}
Let $r$ be a vertex of a semi-2-tree $G$ that is not in any skeleton of $G$.
\mbox{If the root $r$ is} in two fans of $G$, centered on $x$ and $y$, with edge $e=xy$,
then $\pi_t(G,r)=$ $\pi_t(G-e,r)$.\pf
\end{cor}

Since $G-e$ is a semi-2-tree with cut-vertex root $r$, the value
of $\pi_t(G-e,r)$ is computed by Corollary \ref{tPebbSimplicialTree}.

\begin{thm}\label{BestRoot3}
Let $r$ be a vertex of a semi-2-tree $G$ that is not in any skeleton of $G$.
Suppose that the root $r$ is in two fans of $G$, centered on $x$ and $y$, with edge $e=xy$, labeled so that $\ecc(x)\ge\ecc(y)$, and let $r^*$ be a simplicial vertex of $G$ with $\ecc(r^*)=\diam(G)$.
Then $\pi_t(G,r)<\pi_t(G,r^*)$.
\end{thm}

\Pf
We note that $r$ is a cut vertex of $G-e$, and so $b(G-e)=b(G)+1$.
Also, the skeleton $T^\pr$ of $G-e$ has one more edge than does $T$; in fact, $T^\pr$ can be seen as the tree obtained from $T$ by subdividing the edge $e$ with the vertex $r$.
Furthermore, $r\not\in N(r^*)$, which means that $\ecc_{T^\pr}(r)\le\diam(T^\pr)-2$.
As in Lemma \ref{l:tree3}, define $a=\ecc_{T-y}(x)$, and set $d=\diam(T)$.
Then, by Corollaries \ref{tPebbSimplicialTree} and \ref{twofans}, and Lemmas \ref{TreeLemma} and \ref{l:tree3}, we have
\begin{align*}
\pi_t(G,r)
& = \pi_t(G-e,r)\\
& = p_t(G-e,r)\\
& = \pi_t(T^\pr,r)+(n(G-e)-1)+b(G-e)-2e(T^\pr)\\ 
& = \pi_t(T^\pr,r)+(n(G)-1)+b(G)+1-2e(T)-2\\
& = \pi_t(T,x)+2^{a}+(n(G)-1)+b(G)+1-2e(T)-2\\
& \le \pi_t(T,x)+2^{d-\ecc(x)}+(n(G)-1)+b(G)-2e(T)-1\ .\\
\end{align*}
If $\ecc(x)\ge 3$ then we write
\begin{align*}
\pi_t(G,r)
& \le \pi_t(T,x)+2^{d-\ecc(x)}+(n(G)-1)+b(G)-2e(T)-1\\
& \le \pi_t(T,r^*)-(2^{d-2}-2^{d-3})+(n(G)-1)+b(G)-2e(T)-1\\
& < \pi_t(G,r^*)\ .
\end{align*}
Otherwise we have $\ecc(x)\le 2$ and so, with $h^\pr$ defined as in Lemma \ref{TreeLemma}, we find that
\begin{align*}
\pi_t(G,r)
& \le \pi_t(T,x)+2^{d-\ecc(x)}+(n(G)-1)+b(G)-2e(T)-1\\
& \le \pi_t(T,r^*) - t(2^d-2^{\ecc(x)}) + 2^{h^\pr} - 1 + 2^{d-\ecc(x)} +(n(G)-1)+b(G)-2e(T)-1\\
& = \pi_t(G,r^*) - t(2^d-2^{\ecc(x)}) + 2^{h^\pr} - 1 + 2^{d-\ecc(x)}\\
& \le \pi_t(G,r^*) - (2^d - 2^{d-2} - 2^{\lfloor d/2\rfloor} + 1)\\
& < \pi_t(G,r^*)\ ,
\end{align*}
since $d\ge 3$.
\pf

We pause to develop some notation that will be used in Corollary \ref{onefanold}.
Suppose that $r$ is a fan vertex of $G$, in a unique fan $F$ centered on $x$.
Denote by $H_1$ and $H_2$ the two components of $G-\{r,x\}$, and by
$G_i$ the subgraph of $G$ induced by $V(H_i)\cup\{r,x\}$.
Let $V_i$ be the vertices of $F^\pr\cap H_i$ that are not in any other fan.
Define $\Gpi=G_i-V_i$.
Finally, let the subscripts be labeled either so that $V_2$ is empty or so that neither $V_1$ nor $V_2$ is empty and $\ecc_{G_1}(r)\ge\ecc_{G_2}(r)$.

\begin{figure}
\begin{center}

\begin{tikzpicture}[scale=1.2]
\tikzstyle{every node}=[draw,circle,fill=black,minimum size=2pt,inner sep=3pt]
\draw (-5.0,0.0) node (a) {};
\path (-4.0,0.0) node (b) {};
\path (-3.0,0.0) node (c) {};
\path (-2.85,0.6) node (d) {};
\path (-2.35,1.0) node (r) [label=above: $r$] {};
\path (-2.0,0.0) node (x) [label=below: $x$] {};
\path (-1.65,1.0) node (s) {};
\path (-1.15,0.6) node (t) {};
\path (-1.0,0.0) node (u) {};
\path (0.0,0.0) node (v) {};
\draw (a)
  -- (b)
  -- (c)
  -- (d)
  -- (r)
  -- (s)
  -- (t)
  -- (u)
  -- (v);
\draw (x) -- (c);
\draw (x) -- (d);
\draw (x) -- (r);
\draw (x) -- (s);
\draw (x) -- (t);
\draw (x) -- (u);

\draw (1.0,0.0) node (a) {};
\path (2.0,0.0) node (b) {};
\path (3.0,0.0) node (c) {};
\path (3.15,0.6) node (d) {};
\path (3.65,1.0) node (r) [label=above: $r$] {};
\path (4.0,0.0) node (x) [label=below: $x$] {};
\draw (a)
  -- (b)
  -- (c)
  -- (d)
  -- (r);
\draw (x) -- (c);
\draw (x) -- (d);
\draw (x) -- (r);

\path (4.65,1.0) node (r) [label=above: $r$] {};
\path (5.0,0.0) node (x) [label=below: $x$] {};
\path (5.35,1.0) node (s) {};
\path (5.85,0.6) node (t) {};
\path (6.0,0.0) node (u) {};
\path (7.0,0.0) node (v) {};
\draw (r)
  -- (s)
  -- (t)
  -- (u)
  -- (v);
\draw (x) -- (r);
\draw (x) -- (s);
\draw (x) -- (t);
\draw (x) -- (u);

\end{tikzpicture}
\end{center}

\caption{A semi-2-tree $G$ (left), split into $G_1$ (center) and $G_2$ (right).} \label{SplitEx}
\end{figure}
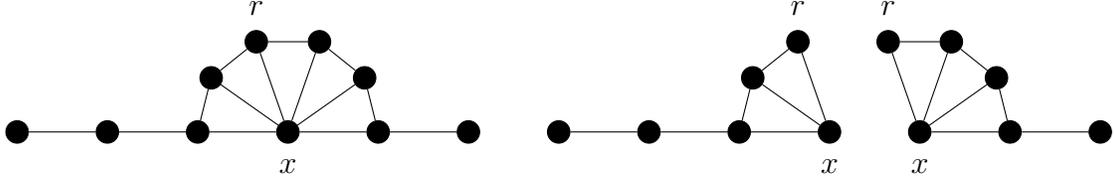

We note that $G_i$ is a semi-2-tree except in the case that the block of $G_i$ containing $r$, $x$ and their unique common neighbor $y$ is a $K_3$ (as in Figure \ref{SmallNbrRem}), because $K_3$ is not a 2-path.
Observe that this happens if and only if $V_i=\emptyset$ and $y$ is a cut vertex of $G$, and that in such a case $G_i-x$ is a semi-2-tree.
Moreover, by the Neighbor Removal Lemma \ref{NeighborRemoval} with $A=\{x\}$, $\pi(G_i,r)=\pi(G_i-x,r)+1$.

\begin{clm}\label{splitformula}
Let $G$ be a semi-2-tree and suppose that $r$ is a fan vertex of $G$, in a unique fan $F$ centered on $x$.
Define $G_i$ ($i\in\{1,2\}$) as above, having $n_i$ vertices.
Define $T_i(r)$ to be the skeleton of $G_i$ when $G_i$ is a semi-2-tree and of $G_i-x$ when $G_i$ is not a semi-2-tree.
Define $T_i(x)$ to be the skeleton of $G_i-V_i-r$.
Then for each $v\in\{r,x\}$ we have $\pi(G_i,v)=\pi(T_i(v),v)+(n_i-1)+b(G_i)-2e(T_i(v))$.
\end{clm}

\Pf
For $v=r$ the result is true by Corollary \ref{tPebbSimplicialTree} when $G_i$ is a semi-2-tree, because $r$ is simplicial.
If $G_i$ is not a semi-2-tree then $|V_i|=0$ and $\{r,x,y\}$ is a $K_3$ block, where $y$ is the common neighbor of $r$ and $x$.
In this case $G_i-x$ is a semi-2-tree, and so $\pi(G_i,r)=\pi(G_i-x,r)+1$ by Lemma \ref{NeighborRemoval}.
This equals $\pi(T_i(r),r)+((n_i-1)-1)+b(G_i-x)-2e(T_i(r))+1
=\pi(T_i(r),r)+(n_i-1)+b(G_i)-2e(T_i(r))$.

For $v=x$ we have that $x$ is a simplicial vertex of the semi-2-tree $G_i-V_i-r$, and so by Lemma \ref{NeighborRemoval} we obtain that
$\pi(G_i,x)
=\pi(G_i-V_i-r,x)+|V_i|+1
=\pi(T_i(x),x)+((n_i-|V_i|-1)-1)+b(G_i-V_i-r)-2e(T_i(x))+|V_i|+1
=\pi(T_i(x),x)+(n_i-1)+b(G_i)-$  $2e(T_i(x))$. \pf

\begin{clm}\label{splitdifference}
Under the same hypotheses as in Claim \ref{splitformula} we have
$\pi(G_1,r)+\pi(G_2,x)\ge\pi(G_1,x)+\pi(G_2,r)$.
\end{clm}

\Pf
Because of the cancellation of common terms, we have\\

$\left[\pi(G_1,r)+\pi(G_2,x)\right]-\left[\pi(G_1,x)+\pi(G_2,r)\right]$
\begin{align*}
& = \left[\pi(T_1(r),r)+\pi(T_2(x),x)\right]-\left[\pi(T_1(x),x)+\pi(T_2(r),r)\right]\\
& = \left[\pi(T_1(r),r)-\pi(T_1(x),x)\right]-\left[\pi(T_2(r),r)-\pi(T_2(x),x)\right]\ .
\end{align*}
Because of the cancellation of common branches, this equals
\begin{align}
\label{sumdiffs}
& \left[2^{\ecc_{G_1}(r)}-2^{\ecc_{G_1}(x)}\right]-\left[2^{\ecc_{G_2}(r)}-2^{\ecc_{G_2}(x)}\right].
\end{align}
We note that $\ecc_{G_i}(x)\le\ecc_{G_i}(r)\le\ecc_{G_i}(x)+1$ for each $i$, with $\ecc_{G_i}(x)=\ecc_{G_i}(r)$ precisely when $V_i=\emptyset$.
Thus, the choice of labeling ensures that (\ref{sumdiffs}) is non-negative.
\pf

\begin{cor}\label{onefanold}
Let $(G,r)$ be a rooted semi-2-tree with $r$ not in any skeleton of $G$.
If $r$ is in a unique fan, centered on $x$, then (using the notation defined above)
$\pi(G,r)=\pi(G_1,r)+\pi(G_2,x)-2$.
\end{cor}


\Pf
The lower bound is argued as follows.
Let $C_1$ be an extremal configuration for $r$ on $G_1$,
$C_2$ be an extremal configuration for $x$ on $G_2-r$
(which is defined by using $A=V_2$ in the Neighbor Removal Lemma \ref{NeighborRemoval}), and define the configuration $C=C_1+C_2$.

Now $|C_1|=\pi(G_1,r)-1$ and $|C_2|=\pi(G_2-r,x)-1=\pi(G_2,x)-2$ (by Lemma \ref{NeighborRemoval} with $A=\{r\}$), and thus $|C|=|C_1|+|C_2|=\pi(G_1,r)+\pi(G_2,x)-3$.
Furthermore, we claim that $C$ is $r$-unsolvable.
Indeed, $C_1$ cannot solve $r$ by itself and cannot receive another
pebble from $C_2$ through $x$, and $C_2$ (without its pebble already
on $r$) cannot solve $r$ by itself (any step to $r$ can be replaced
by a step to $x$, which would be a contradiction).

For the upper bound, assume that $|C|=\pi(G_1,r)+\pi(G_2,x)-2$.
Let $i\in\{1,2\}$ and $j=3-i$.
Define $C_i$ to be the restriction of $C$ to $G_i$.
If $|C_i|\ge \pi(G_i,r)$ then $C_i$ can solve $r$, so we assume otherwise.
Then $|C_j|=|C|-|C_i|+C(x)
\ge \left[\pi(G_1,r)+\pi(G_2,x)-2\right]-\left[\pi(G_i,r)-1\right]+C(x)
\ge \pi(G_j,x)-1+C(x)$.
Indeed, this follows trivially for $j=2$, and from Claim \ref{splitdifference} for $j=1$.
If $C(x)\ge 2$ then we can move a pebble to $r$.
If $C(x)=1$ then, since we may assume that $C(r)=0$,
we have $|C(G_j-r-x)|\ge \pi(G_j,x)-1=\pi(G_j-r,x)$ by Lemma \ref{NeighborRemoval},
and so we can move a second pebble to $x$ and then one to $r$.
Hence we will assume that $C(x)=0$.

If $|C(V_i)|\ge |V_i|+2$ then $V_i$ either has a huge vertex or two big vertices.
in which case it can solve $r$ through $x$, or it has a big vertex with a
path of all ones to $r$, which also solves $r$.
Hence we assume that each $|C(V_i)|\le |V_i|+1$.
Thus we have that
\begin{align*}
|C(\Gpi)| 
 & = |C_i| - |C(V_i)|\\
 & \ge \pi(G_i,x) - |V_i| - 2\\
 & = \pi(\Gpi,x) - 1\ ,
\end{align*}
for each $i$.

Also, if some $|C(\Gpi)|\ge\pi(\Gpi,x)$ then we could place a pebble on $x$.
This implies that $|C(V_j)|\le |V_j|$ since a big vertex in $V_j$ could
place a second pebble on $x$, and then one on $r$.
Then we would have
\begin{align*}
|C(\Gpj)| 
 & = |C_j| - |C(V_j)|\\
 & \ge \pi(G_j,x) - |V_j| - 1\\
 & = \pi(\Gpj,x)\ ,
\end{align*}
so that we could place a second pebble on $x$ and solve $r$.
Thus we must have $|C(\Gpi)|=\pi(\Gpi,x)-1$ for each $i$.

Finally we see that
\begin{align*}
|V_1|+|V_2|+2
 & \ge |C(F)|\\
 & = |C| - |C(\Gpo)| - |C(\Gpt)|\\
 & = [\pi(G_1,r)-\pi(\Gpo,x)] + [\pi(G_2,x)-\pi(\Gpt,x)]\\
 & = [\pi(G_1,r)-\pi(G_1,x)] + [\pi(G_1,x)-\pi(\Gpo,x)] + [|V_2|+1]\\
 & = [2^{\ecc_{G_1}(r)}-2^{\ecc_{G_1}(x)}]+[|V_1|+1] + [|V_2|+1]\ ,
\end{align*}
which means that
$2^{\ecc_{G_1}(x)}\ge 2^{\ecc_{G_1}(r)}$, and hence $2^{\ecc_{G_1}(x)}=2^{\ecc_{G_1}(r)}$.
That is, $V_1=\emptyset$, which implies by our labeling that $V_2=\emptyset$.
Define $x^-$ and $x^+$ to be the common neighbors of $r$ and $x$.
Then in the skeleton of $G$ we can replace the path $x^-xx^+$ by the path $x^-rx^+$ to obtain
a new skeleton containing $r$, which is a contradiction, completing the proof.
\pf


Notice that the previous corollary allows one to calculate the pebbling number for $r$.
In fact, one can use Corollaries \ref{tPebbSimplicialTree} and \ref{spinepebb} to calculate $\pi(G_1,r)$ and $\pi(G_2,x)$, respectively.

As with Corollary \ref{tCheap}, the following is a simple consequence of 
Lemma \ref{Cheap}.

\begin{cor}\label{fanCheap}
Let $(G,r)$ be a rooted semi-2-tree with $r$ not in any skeleton of $G$. If $C$ is a configuration of size at least $\pi(G,r)+(t-1)2^{\ecc(r)}$ then $C$ has $t$ distinct cheap $r$-solutions.
\pf
\end{cor}

Similarly, Corollaries \ref{twofans}, \ref{onefanold}, and \ref{fanCheap}
yield the following result.

\begin{cor}\label{fantpebb}
Let $(G,r)$ be a rooted semi-2-tree with $r$ not in any skeleton of $G$.
Then $\pi_t(G,r)=\pi(G,r)+(t-1)2^{\ecc(r)}$.
\pf
\end{cor}

\begin{thm}\label{BestRoot4}
Let $(G,r)$ be a rooted semi-2-tree with $r$ not in any skeleton of $G$, and let $r^*$ be a simplicial vertex of $G$ with $\ecc(r^*)=\diam(G)$.
Then $\pi_t(G,r) \le \pi_t(G,r^*)$, with equality if and only if $\ecc(r)=\diam(G)$.
\end{thm}

\Pf
We prove that $\pi(G,r) \le\pi(G,r^*)$; then $\pi_t(G,r) = \pi(G,r) + (t-1)2^{\ecc(r)}  \le \pi(G,r^*) + (t-1)2^{\ecc(r^*)} = \pi_t(G,r^*)$ will follow.

First we analyze the case in which $\ecc(r)=\ecc(r^*)$.
Define $x$ to be the center of the fan containing $r$.
Then we can suppose that $x$ is in first (longest) path $P^*$ in the maximum path partition of $T$ with root $ r^*$.
If $s^*$ is the other endpoint of $P^*$ then $x$ is adjacent to $s^*$.
Hence $\ecc_{G_2}(x)=1$ and so $\pi(G_2,x)=n(G_2)=|V_2|+3$.
Thus
\begin{align*}
\pi(G,r)
& = \pi(G_1,r)+\pi(G_2,x)-2\\
& = \pi(G_1,r)+| V_2|+3-2\\
& = \pi(G_1,r)+|V_2|+1 \\
\end{align*}
Also,
\begin{align*}
\pi(G,r^*)
& = \pi(T,r^*) + (n(G)-1) + b(G) - 2e(T)\\
& = \pi(T_1,r)+ (n(G_1)+ |V_2|+1 -1) + b(G_1) - 2e(T_1) \\
& = \pi(G_1,r)+ |V_2| +1\ .\\
\end{align*}

Henceforth we will assume that $\ecc(r)<\ecc(r^*)$. 
In this case we will make use of Lemma \ref{TreeLemma} and Corollaries \ref{onefanold} and \ref{fantpebb}.
We will also use the facts that $n(G_1)+n(G_2)=n(G)+2$, $b(G_1)+b(G_2)=b(G)+1$, and $e(T_1)+e(T_2)=e(T)+\epsilon$, where $\epsilon$, depends on some cases.

Notice that $T_2$ does not include the edge $xr$ because $x$ is the root, while $T_1$ does include the edge $xr$ unless $V_1=\emptyset$ and $r$ and $x$ have a common neighbor in $T\cap V_1$.
Hence $\epsilon=1$ except in this latter case, in which $\epsilon=0$; that is, $\epsilon=1-|N(r)\cap N(x)\cap T\cap V_1|$.
Now

\begin{align*}
\pi(G,r)
& = \pi(G_1,r)+\pi(G_2,x)-2\\
& = \pi(T_1,r) + (n(G_1)-1)+b(G_1)-2e(T_1)\\
& \qquad + \pi(T_2,x) + (n(G_2)-1)+b(G_2)-2e(T_2) - 2\\
& = \pi(T_1,r) + \pi(T_2,x) + (n(G)-1) + b(G) - 2e(T) - 2\epsilon\ .\\
\end{align*}

Analogous to the proof of Lemma \ref{TreeLemma}, let $P^*$ be a path $v_0v_1\cdots v_d$ with $v_0=r^*$ and $v_d=s^*$, labeled so that $\dist(x,s^*)\le\dist(x,r^*)=\ecc(x)$.
Denote by $P$ the path from $r^*$ to $x$, and set $P^*\cap P=v_0\cdots v_{h^\pr}$.
Define $\overline{h}=d-h^\pr$.
let $\cP^*$ be a maximum path partition of $T$ with root $r^*$.
Define $P^*_0=P^*$, $P^*_1$, $\ldots$, $P^*_k$ to be the sequence of paths of $\cP^*$ that are used sequentially in $P$, and set $d^*_i=\len(P^*_i)$ for each $0\le i\le k$ (so $d^*_0=d$).
Next define $P_i^\pr=P\cap P^*_i$, with $h^\pr_i=\len(P_i^\pr)$ and $\overline{h}_i=d^*_i-h^\pr_i$ (so $h^\pr_0=h^\pr$ and $\overline{h}_0=\overline{h}$).
Note that $1\le h^\pr_i\le d^*_i\le d/2$.

Suppose that $k>0$.
Then, since $\ecc(r)<\ecc(r^*)$, we have
\begin{align*}
\pi(T_1,r) + \pi(T_2,x)
& \le \pi(T_1,r^*) - 2^{d-2} + \pi(T_2,x)\\
& = \pi(T,r^*) + (2^{h^\pr_k+1}-1) + (2^{\overline{h}_k}-1) - (2^{d^*_k}-1) - 2^{d-2}\\
& \le \pi(T,r^*) + 1 - 2^{d-2}\\
& < \pi(T,r^*)\ ,
\end{align*}
because $d\ge 3$ and $2^{a+1}+2^b-2^{a+b}\le 2$ for all $a, b\ge 1$.

Suppose instead that $k=0$.
Define $Q$ to be the longest path in a maximum path partition of $T_1$ (choosing the partition to contain $P$, if possible).
If $Q=P$ then we have
\begin{align*}
\pi(T_1,r) + \pi(T_2,x)
& \le \pi(T_1,r^*) - 2^{(h^\pr_0+1)-2} + \pi(T_2,x)\\
& = \pi(T,r^*) - (2^d-1) + (2^{h^\pr_0+1}-1) + (2^{\overline{h}_0}-1) - 2^{h^\pr_0-1}\\
& < \pi(T,r^*) - 2^d - 2^{h^\pr_0+1} + 2^{\overline{h}_0} - 2^{h^\pr_0-1}\\
& < \pi(T,r^*) - 2^d + 2^{\overline{h}_0}\\
& \le \pi(T,r^*) - 2^d + 2^{d-1}\\
& < \pi(T,r^*)\ .
\end{align*}
Otherwise, when $Q\not=P$ we define $Q_0=Q\cap P$ and $Q_1=Q-Q_0$, having lengths $q_0$ and $q_1$, respectively, and set $\hat{h}_0=h^\pr_0-q_0$.
Here we will use that $q_0+q_1<d$, $q_1>\hat{h}_0$, and $\overline{h}_0\le d/2\le q_0+\hat{h}_0\le q_0+q_1$.
Hence
\begin{align*}
\pi(T_1,r) + \pi(T_2,x)
& \le \pi(T_1,r^*) - 2^{q_0+q_1-2} + \pi(T_2,x)\\
& = \pi(T,r^*) - (2^d-1) - (2^{q_1}-1) + (2^{q_0+q_1}-1) + (2^{\hat{h}_0+1}-1)\\ 
& \qquad + (2^{\overline{h}_0}-1) - 2^{q_0+q_1-2}\\
& < \pi(T,r^*) - 2^d  - 2^{q_1} + 2^{q_0+q_1} + 2^{\hat{h}_0+1} + 2^{q_0+q_1} - 2^{q_0+q_1-2}\\
& < \pi(T,r^*) - (2^d  - 2^{q_0+q_1+1}) - (2^{q_1} - 2^{\hat{h}_0+1})\\
& \le \pi(T,r^*)\ .
\end{align*}

In all cases, then, we see that 
\begin{align*}
\pi(G,r)
& = \pi(T_1,r) + \pi(T_2,x) + (n(G)-1) + b(G) - 2e(T) - 2\epsilon\ .\\
& < \pi(T,r^*) + (n(G)-1) + b(G) - 2e(T)\\
& = \pi(G,r^*)\ ,
\end{align*}
and the result follows.
\pf

\begin{thm}\label{S2Tformula}
If $G$ is a semi-2-tree then $\pi_t(G)=\pi_t(G,r^*)$, where $r^*$ is a simplicial vertex with $\ecc(r^*)=\diam(G)$.
\end{thm}

\Pf
Use Theorems \ref{BestRoot1}, \ref{BestRoot2}, \ref{BestRoot3}, and \ref{BestRoot4}.
\pf

\begin{thm}\label{S2Tcomplexity}
If $G$ is a semi-2-tree then $\pi_t(G)$ can be computed in
linear time.
\end{thm}

\Pf
A breadth-first search from any simplicial vertex finds $r^*$, a simplicial vertex with $\ecc(r^*)=\diam(G)$.
Indeed, this is true for trees, and the result extends to semi-2-trees as follows.
Let $T$ be the skeleton of $G$ and let $\mathbf{A}$ be a breadth-first search algorithm on $G$.
Then $\mathbf{A}$ is also a breadth-first search algorithm on $T$ and so finds a simplicial vertex $r$ with $\ecc_T(r)=\diam(T)$.
Because $T$ is a geodesic tree spanning all of the simplicial vertices of $G$, we have $\ecc_G(r)=\ecc_T(r)$ and $\diam(G)=\diam(T)$, and so $r^*=r$.

At this point, we do not yet know $T$.
However, we realize that $T$ can be constructed during $\mathbf{A}$ because it is a geodesic tree spanning all of the simplicial vertices of $G$.
Once we have $T$ we can remove its cut-vertices $S$ (those having degree bigger than 2) to reveal $b$, which equals the number of components of $T-S$.

Then $\pi_t(T,r)$ can be computed in linear time, according to Theorem 3 of \cite{BCCMW}.
\pf

%
%
\section{Remarks}\label{Remarks}

The obvious pressing question is how to extend this work to $2$-trees.
The {\it pyramid} is the graph on 6 vertices formed by adjoining a
2-simplicial vertex onto each of the three sides of a triangle.  The
pyramid is the key structure that forms the basis in the Class 0 
characterization of diameter two graphs found in \cite{ClHoHu} and
is what causes the extra 1 in their pebbling numbers --- the configuration
with 3 pebbles at two of the simplicial vertices cannot reach the third.
The pyramid is also the smallest example of a 2-tree that is not a 
semi-2-tree, and it hints at the complexity that can ensue in a more
general 2-tree.

Another natural question in the direction of this research program
regards other simple examples of chordal graphs, such as interval
graphs.  It would seem that tackling $k$-paths is a necessary
investigation toward approaching interval graphs.  One interesting
thing about the 2-path pebbling number is that both of the standard
lower bounds of $n(G)$ and $2^{\diam(G)}$ for general graphs $G$
appear in its formula.  This is encouraging in light of the manner
in which the size of $k$ can determine which of those two terms is
dominant.

It appears that parameters such as pathwidth and treewidth may figure
prominantly in the determination of pebbling numbers of general graphs.  
Other authors have made similar remarks, for example in \cite{ChaGod}.  Thus
considering these classes of graphs seems the most productive direction of
research.

Our final thought points to the many lemmas developed in this paper
that should be of very general use, including the Cheap Lemma
(\ref{Cheap}) and the four Removal lemmas: Junior (\ref{JuniorRemoval}),
Wart (\ref{WartRemoval}), Edge (\ref{EdgeRemoval}), and Neighbor
(\ref{NeighborRemoval}).  We anticipate their ability to simplify
the analysis of many future problems.

\section{Acknowledgment}

We wish to thank the referees for excellent feedback that
improved the exposition in this article greatly.

%
%

\end{document}